\documentclass[12pt]{amsart}%
\usepackage{graphicx}
\usepackage{amscd}
\usepackage{amsmath}
\usepackage{amsfonts}
\usepackage{amssymb}%
\providecommand{\U}[1]{\protect\rule{.1in}{.1in}}
\theoremstyle{plain}

\theoremstyle{definition}

\newtheorem{remark}{Remark}[section]

\numberwithin{equation}{section}
\numberwithin{theorem}{section}
\begin{document}
\title[Orthogonal polynomials]{Note on the $X_{1}$-\textbf{Jacobi} orthogonal polynomials.}
\author{W.N. Everitt}
\address{W.N. Everitt, School of Mathematics, University of Birmingham, Edgbaston,
Birmingham B15 2TT, England, UK}
\email{w.n.everitt@bham.ac.uk}
\date{21 November 2008 (File C:%
$\backslash$%
Swp55%
$\backslash$%
Docs%
$\backslash$%
milson23.tex)}

\begin{abstract}
This note supplements the results in the paper on $X_{1}$-\textbf{Jacobi}
orthogonal polynomials, written by David G\'{o}mez-Ullate, Niky Kamran and
Robert Milson.

\end{abstract}
\subjclass[2000]{ Primary; 34B24; 34L05, 33C45: Secondary; 05E35, 34B30.}
\keywords{Sturm-Liouville theory, orthogonal polynomials.}
\maketitle

\section{Introduction\label{sec1}}

This note reports on, the $X_{1}$-\textbf{Jacobi polynomials}, one of the two
new sets of orthogonal polynomials considered in the papers \cite{G-UKM} and
\cite{G-UKM1}, written by David G\'{o}mez-Ullate, Niky Kamran and Robert
Milson. The other set is named the $X_{1}$-\textbf{Laguerre polynomials} and
is discussed, in similar terms, in the note \cite{WNE1}.

These two papers are remarkable and invite comments on the results therein
which have yielded new examples of Sturm-Liouville differential equations and
their associated differential operators.

The two sets of these orthogonal polynomials are distinguished by:

\begin{itemize}
\item[$(i)$] Each set of polynomials is of the form $\{P_{n}(x):x\in
\mathbb{R}\ $and$\ n\in\mathbb{N}\equiv\{1,2.3.\ldots\}\}$ with$\ \deg
(P_{n})=n;$ that is there is no polynomial of degree $0.$

\item[$(ii)$] Each set is orthogonal and complete in a weighted Hilbert
function space.

\item[$(iii)$] Each set is generated as a set of eigenvectors from a
self-adjoint Sturm-Liouville differential operator.
\end{itemize}

\section{$X_{1}$-\textbf{Jacobi Polynomials\label{sec2}}}

\subsection{Parameters\label{sec2.1}}

These polynomials and associated differential equations are detailed in
\cite[Section 2]{G-UKM}

There are two real-valued parameters $\alpha,\beta$ involved in the
differential equation for these orthogonal polynomials. These parameters have
to satisfy, see \cite[Section 2.1]{G-UKM},:%
\begin{equation}
\alpha,\beta\in(-1,\infty)\hspace{0.5in}\alpha\neq\beta\hspace{0.5in}%
\mathrm{sgn}(\alpha)=\mathrm{sgn}(\beta) \label{eq2.0}%
\end{equation}

These conditions give two essentially different cases to consider%
\begin{equation}
\left\{
\begin{array}
[c]{cc}%
\text{\textbf{Case}}\ 1 & \beta>\alpha>0\\
\text{\textbf{Case}}\ 2 & -1<\beta<\alpha<0
\end{array}
\right.  \label{eq2.1}%
\end{equation}
The real parameters $a,b,c,$ see \cite[Section 2, (5a) and (5b)]{G-UKM}, are
defined by%
\begin{equation}
a:=\tfrac{1}{2}(\beta-\alpha),\ b:=\dfrac{\beta+\alpha}{\beta-\alpha
},\ c:=b+a^{-1}; \label{eq2.2}%
\end{equation}
it may be shown, using (\ref{eq2.0}), that $a,b,c,$ satisfy%
\begin{equation}
\left\{
\begin{array}
[c]{cc}%
\text{\textbf{Case}}\ 1 & a>0,\ b>1,\ c>1\\
\text{\textbf{Case}}\ 2 & a<0,\ b>1,\ c<-1.
\end{array}
\right.  \label{eq2.3}%
\end{equation}

\subsection{Differential expression\label{sec2.2}}

In \cite[Section 2, (19) and (20a)]{G-UKM} the second-order linear
differential equation concerned is given as%
\begin{equation}
(x^{2}-1)y^{\prime\prime}(x)+2a\left(  \dfrac{1-bx)}{b-x}\right)  \left[
(x-c)y^{\prime}(x)-y(x)\right]  =\lambda y(x)\ \text{for all}\ x\in(-1,1),
\label{eq2.4}%
\end{equation}
where the parameter $\lambda\in\mathbb{C}$ plays the role of a spectral
parameter for the differential operators defined below.

This equation (\ref{eq2.4}) is not a Sturm-Liouville differential equation;
such equations take the form, in this case taking the interval to be $(-1,1),$%
\begin{equation}
-(p(x)y^{\prime}(x))^{\prime}+q(x)y(x)=\lambda w(x)y(x)\ \text{for all}%
\ x\in(-1,1). \label{eq2.5}%
\end{equation}
with $\lambda\in\mathbb{C}.$ Here $p,q,w:(-1,1)\rightarrow\mathbb{R}$ and
satisfy the minimal conditions, \cite[Section 3]{WNE},%
\begin{equation}%
\begin{array}
[c]{ll}%
(i) & p^{-1},q,w\in L_{\text{loc}}^{1}(-1,1)\\
(ii) & w(x)>0\ \ \text{for almost all}\ x\in(-1,1).
\end{array}
\label{eq2.6}%
\end{equation}

The equation (\ref{eq2.4}) can be transformed into Sturm-Liouville form on
multiplication by the weight $\hat{W}_{\alpha,\beta}$ given in \cite[Section
2, (11)]{G-UKM}, but here denoted by $w_{\alpha,\beta}$ and defined by%
\begin{equation}
w_{\alpha,\beta}(x):=\dfrac{(1-x)^{\alpha}(1+x)^{\beta}}{(x-b)^{2}}\ \text{for
all}\ x\in(-1,1). \label{eq2.7}%
\end{equation}
If now the equation (\ref{eq2.4}) is multiplied by $w_{\alpha,\beta}$ then the
Sturm-Liouville form (\ref{eq2.5}) is satisfied with%
\begin{equation}
p_{\alpha,\beta}(x):=\dfrac{(1-x)^{\alpha+1}(1+x)^{\beta+1}}{(x-b)^{2}}\text{
for all}\ x\in(-1,1) \label{eq2.8}%
\end{equation}
and%
\begin{equation}
q_{\alpha,\beta}(x):=2a\left(  \dfrac{1-bx)}{b-x}\right)  (x-c)\dfrac
{(1-x)^{\alpha}(1+x)^{\beta}}{(x-b)^{2}}\ \text{for all}\ x\in(-1,1).
\label{eq2.9}%
\end{equation}

In passing we note that, for the chosen values of $\alpha,\beta$ in
(\ref{eq2.1}), the Sturm-Liouville equation, defining the differential
expression $M_{\alpha,\beta},$%
\begin{equation}
M_{\alpha,\beta}[y](x):=-(p_{\alpha,\beta}(x)y^{\prime}(x))^{\prime}%
+q_{\alpha,\beta}(x)y(x)=\lambda w_{\alpha,\beta}(x)y(x)\ \text{for all}%
\ x\in(-1,1). \label{eq2.10}%
\end{equation}
is regular at all points of the open interval $(-1,1)$.

For \textbf{Case} $1$ of (\ref{eq2.3}) the differential equation
(\ref{eq2.10}) is singular at both endpoints $\pm1$ since%
\[
\int_{0}^{1}\dfrac{1~}{p_{\alpha,\beta}(x)}~dx=\int_{-1}^{0}\dfrac
{1~}{p_{\alpha,\beta}(x)}~dx=+\infty.
\]

For \textbf{Case} $2$ of (\ref{eq2.2}) the differential equation
(\ref{eq2.10}) is regular at both endpoints $\pm1$ since%
\[
\int_{0}^{1}\dfrac{1~}{p_{\alpha,\beta}(x)}~dx\,<+\infty\ \text{and}%
\ \int_{-1}^{0}\dfrac{1~}{p_{\alpha,\beta}(x)}~dx<+\infty.
\]

The symplectic form for $M_{\alpha,\beta}$ is defined by, for all
$k\in(0,\infty)$ and for all $f,g\in D(M_{\alpha,\beta}),$%
\begin{equation}
\lbrack f,g]_{\alpha,\beta}(x):=f(x)(p_{\alpha,\beta}\overline{g}^{\prime
})(x)-(p_{\alpha,\beta}f^{\prime})(x)\overline{g}(x)\ \text{for all}%
\ x\in(-1,1). \label{eq2.11}%
\end{equation}

\section{Differential operators\label{sec3}}

The Sturm-Liouville differential expression $M_{\alpha,\beta}$ defines
differential operators in the Hilbert function space%
\[
L^{2}((-1,1);w_{\alpha,\beta}).
\]
In this space the maximal domain of the differential expression $M_{\alpha
,\beta}$ is defined by%
\begin{equation}
\left\{
\begin{array}
[c]{ll}%
(i) & T_{\alpha,\beta}:D(T_{\alpha,\beta})\subset L^{2}((-1,1);w_{\alpha
,\beta})\rightarrow L^{2}((-1,1);w_{\alpha,\beta})\\
(ii) & D(T_{\alpha,\beta}):=\{f\in D(M_{\alpha,\beta}):f,w^{-1}M_{\alpha
,\beta}[f]\in L^{2}((-1,1);w_{\alpha,\beta})\\
(iii) & T_{\alpha,\beta}f:=w^{-1}M_{\alpha,\beta}[f]\ \text{for all}\ f\in
D(T_{\alpha,\beta}).
\end{array}
\right.  \label{eq3.1}%
\end{equation}

All self-adjoint differential operators in $L^{2}((-1,1);w_{\alpha,\beta})$
generated by $M_{\alpha,\beta}$ are given by restrictions of the maximal
operator $T_{\alpha,\beta};$ these restrictions are determined by placing
boundary conditions at the endpoints $-1$ and $+1,$ on the elements of
$D(T_{\alpha,\beta}).$ The number and type of boundary conditions depends upon
the endpoint classification of $M_{\alpha,\beta}$ in $L^{2}((0,\infty
);w_{k});$ see \cite[Section 5]{WNE}.

For the endpoint classification of the differential expression $M_{\alpha
,\beta}$ in $L^{2}((-1,1);w_{\alpha,\beta})$ we have the results, see again
\cite[Section 5]{WNE};

For \textbf{Case} $1$ we have

\begin{itemize}
\item[$(i)$] Endpoint $+1$%
\begin{equation}
\left\{
\begin{array}
[c]{ll}%
\text{limit-circle} & \text{for}\ 0<\alpha<1\\
\text{limit-point} & \text{for}\ \alpha\geq1.\
\end{array}
\right.  \label{eq3.1a}%
\end{equation}

\item[$(ii)$] Endpoint $-1$%
\begin{equation}
\left\{
\begin{array}
[c]{ll}%
\text{limit-circle} & \text{for}\ 0\,<\beta<1\\
\text{limit-point} & \text{for}\ \beta\geq1.\
\end{array}
\right.  \label{eq3.1b}%
\end{equation}

\end{itemize}

The proof of these last results follows using linearly independent solutions
$\varphi_{1}$ and $\varphi_{2}$ of the equation (\ref{eq2.10}) for
$\lambda=0,$ see \cite[Section 2, (6)]{G-UKM},%
\begin{equation}
\varphi_{1}(x):=x-c\ \text{for all}\ x\in(-1,1) \label{eq3.2}%
\end{equation}
and%
\begin{align}
\varphi_{2}(x)  &  :=\varphi_{1}(x)\int_{0}^{x}\dfrac{1}{p_{\alpha,\beta
}(t)\varphi_{1}^{2}(t)}~dt\ \nonumber\\
&  =(x-c)\int_{0}^{x}\dfrac{(t-b)^{2}}{(1-t)^{\alpha+1}(1+t)^{\beta
+1}(t-c)^{2}}\ dt\text{ for all}\ x\in(-1,1). \label{eq3.3}%
\end{align}
Clearly then, recalling the restrictions on the parameters $\alpha
,\beta,a,b,c,$ asymptotic analysis shows that%
\begin{equation}
\varphi_{1}\in L^{2}((-1,1);w_{k})\ \text{for all}\ \alpha,\beta\label{eq3.4}%
\end{equation}
and%
\begin{equation}
\left\{
\begin{array}
[c]{ll}%
\varphi_{2}\in L^{2}((0,1);w_{k}) & \text{for all}\ 0<\alpha<1\\
\varphi_{2}\notin L^{2}((0,1);w_{k}) & \text{for all}\ \alpha\geq1
\end{array}
\right.  \label{eq3.5}%
\end{equation}
~%
\begin{equation}
\left\{
\begin{array}
[c]{ll}%
\varphi_{2}\in L^{2}((-1,0);w_{k}) & \text{for all}\ 0<\beta<1\\
\varphi_{2}\notin L^{2}((-1,0);w_{k}) & \text{for all}\ \beta\geq1.
\end{array}
\right.  \label{eq3.6}%
\end{equation}
These results imply the endpoint conditions stated (\ref{eq3.1a}) and
(\ref{eq3.1b}).

For \textbf{Case} $2$ we have both $\varphi_{1}$ and $\varphi_{2}$ belong to
$L^{2}((-1,1);w_{k})$ for all $\alpha,\beta$ and this observation gives%
\begin{equation}
\text{limit-circle for }-1<\beta<\alpha<0. \label{eq3.6c}%
\end{equation}

To determine the restriction $A_{\alpha,\beta}$ of the maximal operator
$T_{\alpha,\beta}$ which yields the $X_{1}$-\textbf{Jacobi} orthogonal
polynomials as eigenvectors we take the boundary condition function at $\pm1$
to be the solution $\varphi_{1},$ and use the symplectic form (\ref{eq2.11});
here $f$ is any element of $D(T_{\alpha,\beta})$;

\textbf{Case} $1$

\begin{itemize}
\item[$(i)$] Endpoint $+1$%
\begin{equation}
\left\{
\begin{array}
[c]{ll}%
\text{for}\ 0\leq\alpha<1 & \lim_{x\rightarrow+1^{-}}\left[  f,\varphi
_{1}\right]  (x)=0\\
\text{or equivalently} & \lim_{x\rightarrow+1^{-}}(1-x)^{\alpha+1}%
(f(x)-f^{\prime}(x)(x-c))=0
\end{array}
\right.  \label{eq3.6a}%
\end{equation}%
\[%
\begin{array}
[c]{cc}%
\text{for}\ \alpha\geq1 & \text{no boundary condition required.}%
\end{array}
\]

\item[$(ii)$] Endpoint $-1$%
\begin{equation}
\left\{
\begin{array}
[c]{ll}%
\text{for}\ 0<\beta<1 & \lim_{x\rightarrow+1^{-}}\left[  f,\varphi_{1}\right]
(x)=0\\
\text{or equivalently} & \lim_{x\rightarrow-1^{+}}(1-x)^{\beta+1}%
(f(x)-f^{\prime}(x)(x-c))=0.
\end{array}
\right.  \label{eq3.6b}%
\end{equation}%
\[%
\begin{array}
[c]{cc}%
\text{for}\ \beta\geq1 & \text{no boundary condition required.}%
\end{array}
\]

\end{itemize}

\textbf{Case} $2$%
\begin{equation}%
\begin{array}
[c]{ll}%
\text{for}\ -1<\,\beta<\alpha<\,0 & \lim_{x\rightarrow+1^{\pm}}\left[
f,\varphi_{1}\right]  (x)=0.
\end{array}
\label{eq3.6d}%
\end{equation}

The domain $D(A_{\alpha,\beta})$ of the self-adjoint restriction
$A_{\alpha,\beta}$ is then determined by applying the above boundary
conditions in the appropriate cases for the parameters $\alpha,\beta$ to give%
\begin{equation}
A_{\alpha,\beta}f:=w_{\alpha,\beta}^{-1}M_{\alpha,\beta}[f]\ \text{for
all}\ f\in D(A_{\alpha,\beta}). \label{eq3.7}%
\end{equation}

The spectrum and eigenvectors of $A_{\alpha,\beta}$ can be obtained from the
results given in \cite[Section 2]{G-UKM}. The spectrum of $A_{\alpha,\beta}$
contains the sequence $\{\lambda_{n}=n(\alpha+\beta+n):n\in\mathbb{N}_{0}\};$
the eigenvectors are given by $\{\hat{P}_{n}^{(\alpha,\beta)}:n\in
\mathbb{N}_{0}\},$ the $X_{1}$-\textbf{Jacobi} orthogonal polynomials.

\begin{remark}
\label{rem2.1}

\begin{itemize}
\item[$(i)$] The notation $\lambda_{n}=n(\alpha+\beta+n)$ for all
$n\in\mathbb{N}_{0}$ makes good comparison with the eigenvalue notation for
the classical Jacobi polynomials; this sequence depends upon the parameters
$\alpha,\beta.$

\item[$(ii)$] We note that $\hat{P}_{n}^{(\alpha,\beta)}$ is a polynomial of
degree $n+1$ for all $n\in\mathbb{N}_{0}$ and all $\alpha,\beta\ $under consideration.

\item[$(iii)$] Note that when the limit-circle condition holds at $\pm1,$ it
is essential to check that the polynomials $\{\hat{P}_{n}^{(\alpha,\beta)}\}$
all satisfy the boundary conditions at $\pm1$ as required in $(\ref{eq3.6a})$,
$(\ref{eq3.6b})$ and $(\ref{eq3.6d}).$ Thus it is required that%
\[
\lim_{x\rightarrow0+}\left[  \hat{P}_{n}^{(\alpha,\beta)},\varphi_{1}\right]
(x)=0\ \text{for all}\ n\in\mathbb{N}_{0.}%
\]
This result follows since, using $(\ref{eq2.8}),$%
\begin{align*}
\left[  \hat{P}_{n}^{(\alpha,\beta)},\varphi_{1}\right]  (x)  &
=p_{\alpha,\beta}(x)\left[  \hat{P}_{n}^{(\alpha,\beta)}(x)\varphi_{1}%
^{\prime}(x)-\hat{P}_{n}^{(\alpha,\beta)\prime}(x)\varphi_{1}(x)\right] \\
&  =\dfrac{(1-x)^{\alpha+1}(1+x)^{\beta+1}}{(x-b)^{2}}\left[  \hat{P}%
_{n}^{(\alpha,\beta)}(x)-\hat{P}_{n}^{(\alpha,\beta)\prime}(x)(x+k+1)\right]
\\
&  =\mathcal{O}((1-x)^{\alpha+1}(1+x)^{\beta+1})\ \text{as}\ x\rightarrow
+1^{\pm}.
\end{align*}

\end{itemize}
\end{remark}

It is shown in \cite[Section 3, Proposition 3.2]{G-UKM} that the sequence of
polynomials
\[
\left\{  \hat{P}_{n}^{(\alpha,\beta)}:n\in\mathbb{N}_{0}\right\}
\]
is orthogonal and dense in the space $L^{2}((-1,1);w_{k}),$ for all
$\alpha,\beta.$ This result implies that the spectrum of the operator
$A_{\alpha,\beta}$ consists entirely of the sequence of eigenvalues
$\{\lambda_{n}=n(\alpha+\beta+n):n\in\mathbb{N}_{0}\};$ from the spectral
theorem for self-adjoint operators in Hilbert space it follows that no other
point on the real line $\mathbb{R}$ can belong to the spectrum of
$A_{\alpha,\beta}.$

\begin{remark}
\label{rem2.2}It is to be noted that whilst the Hilbert space theory as given
in \cite{WNE} and \cite{MAN} provides a precise definition of the self-adjoint
operator $A_{\alpha,\beta},$ the information about the particular spectral
properties of $A_{\alpha,\beta}$ are to be deduced from the classical analysis
results in \cite{G-UKM}. Without these results it would be very difficult to
deduce the spectral properties of the self-adjoint operator $A_{\alpha,\beta
},$ as defined above, in the Hilbert function space $L^{2}((-1,1);w_{k}).$
\end{remark}


\begin{thebibliography}{9}                                                                                                %


\bibitem {WNE}W.N. Everitt. A catalogue of Sturm-Liouville differential
equations. \emph{Sturm-Liouville Theory, Past and Present}: Pages 271-331.
(Birkh\"{a}user Verlag, Basel: 2005; edited by W.O. Amrein, A.M. Hinz and D.B. Pearson.)

\bibitem {WNE1}W.N. Everitt. Note on the $W_{1}$-\textbf{Laguerre} orthogonal
polynomials. (Submitted to arXiv [math-ph] 21 November 2008).

\bibitem {G-UKM}D. G\'{o}mez-Ullate, N. Kamran and R. Milson. An extended
class of orthogonal polynomials defined by a Sturm-Liouville
problem.(arXiv:080/.3939v1 [math-ph] 24 July 2008).

\bibitem {G-UKM1}D. G\'{o}mez-Ullate, N. Kamran and R. Milson. An extension of
Bochner's problem: exceptional invariant subspaces. (arXiv:0805.3376v2
[math-ph] 24 July 2008).

\bibitem {MAN}M.A. Naimark. \emph{Linear differential operators}: Part
\textbf{II}. (Ungar New York: 1968.)
\end{thebibliography}
\end{document}